\documentclass[twoside,12pt,a4paper]{article}
\usepackage{latexsym,amsfonts,amstext,amsbsy,amsmath,amssymb,amscd,graphicx,epsfig,eucal,bbm}
\usepackage{mathrsfs}
\usepackage{pb-diagram}
\usepackage[all]{xy}
\usepackage{makeidx}
\usepackage{indentfirst}
\usepackage{multicol}
\usepackage{graphicx}

\newtheorem{theorem}{Theorem}

\newtheorem{pr}[theorem]{Proposition}

\date{}

\usepackage{fancyhdr}
\begin{document}

\title{\textbf{{STOCHASTIC FRACTIONAL HP EQUATIONS}}}
\author{Chi\c{s} Oana$^1$, Opri\c{s} Dumitru$^2$\\
West University of Timi\c{s}oara, Romania\\
4 Vasile P\^{a}rvan Blvd., Timi\c{s}oara, 300223, Romania\\
chisoana@yahoo.com\\
opris@math.uvt.ro}

\maketitle

\thispagestyle{fancy}

\textbf{Abstract:} In this paper we established the condition for
a curve to satisfy stochastic fractional HP (Hamilton-Pontryagin)
equations. These equations are described using It\^o integral. We
have also considered the case of stochastic fractional Hamiltonian
equations, for a hyper\-regular Lagrange function. From the
stochastic fractional Hamiltonian equations, Langevin fractional
equations were found and numerical simulations were done.\\


\section{Introduction}

J.M. Bismut was the first one that introduced concepts of
stochastic geometric mechanics, in his work from 1981, when he
defined the notion of "stochastic Hamiltonian system". He showed
that that the stochastic flow of a certain randomly perturbed
Hamiltonian systems on flat spaces extremizes a stochastic action,
and using this property, he proved symplecticity and Noether
theorem for stochastic Hamiltonian systems. Since then, there has
been a need in finding out tools and algorithms for the study of
this kind of systems with uncertainty. Bismut's work was continued
by Lazaro-Cami and Ortega (\cite{lazaro}, \cite{lazaro1}), in the
sense that his work was generalized to manifolds, stochastic
Hamiltonian systems on manifolds extremize a stochastic action on
the space of manifold valued semimartingales, the reduction of
stochastic Hamiltonian system on cotangent bundle of a Lie group,
a counter example for the converse of Bismut's original theorem.

Very important in many science fields is fractional calculus:
fractional derivatives, fractional integrals, of any real or
complex order. Fractional calculus is used when fractional
integration is needed. It is used for studying simple dynamical
systems, but it also describes complex physical systems. For
example, applications of the fractional calculus can be found in
chaotic dynamics, control theory, stochastic modelling, but also
in finance, hydrology, biophysics, physics, astrophysics,
cosmology and so on (\cite{chis}, \cite{el-b1}, \cite{el}). But
some other fields have just started to study problems from
fractional point of view. In great fashion is the study of
fractional problems of the calculus of variations and
Euler-Lagrange type equations. There were found Euler-Lagrange
equations with fractional derivatives, and then Klimek found
Euler-Lagrange equations, but with symmetric fractional
derivatives \cite{klimek}. Most famous fractional integral are
Riemann-Liouville, Caputo, Grunwald-Letnikov and most frequently
used is Riemann-Liouville fractional integral. The study of
Euler-Lagrange fractional equations was continued  by Agrawal
\cite{agrawal1} that described these equations using the left,
respectively right fractional derivatives in the Riemann-Liouville
sense. This fractional calculus has some great problems, such as
presence of non-local fractional differential operators, or the
adjoint fractional operator that describes the dynamics is not the
negative of itself, or mathematical calculus may be very hard
because of the complicated Leibniz rule, or the absence of chain
rule, and so on. After O.P. Agrawal's formulation \cite{agrawal1}
of Euler-Lagrange fractional equations, B\u{a}leanu and Avkar
\cite{baleanu} used them in formulating problems with Lagrangians
linear in velocities. Standard multi-variable variational calculus
has also some limitations. But in \cite{udriste} C. Udri\c{s}te
and D. Opri\c{s} showed that these limitations can by broken using
the multi-linear control theory.

For fractional stochastic integrals of the form
\begin{equation}
X_t=X_0+\int_0^tF(t,s,X_s)ds+\int_0^tG_a(H_t;t,s,X_s)dW_s^a,
\end{equation}
the existence and uniqueness of its solution was discussed by
Pardoux and Protter in their work \cite{pardoux}.

In \cite{el-b1}, it was proved the existence, uniqueness, and
continuity of a fractional stochastic equation of the form
\begin{equation}\label{eq2}
X_t=X_0+I_t^{\beta}F(t,s,X_s)+W_t^{\beta}G^a(H_t;t,s,X_s),
\end{equation}where $0<\beta<1,$ $I_t^\beta F(t,s,X_s).$ The
fractional integral of $F(t,s,X_s)$ is defined by
$$I_t^{\beta}F(t,s,X_s)=\frac{1}{\Gamma(\beta)}\int_0^t\frac{F(t,s,X_s)}{(t-s)^{1-\beta}}ds,$$
and fractional Wiener process $W_t^\beta G_a(H_t;t,s,X_s)$ of
$G_a(H_t;t,s,X_s)$ is
$$W_t^\beta G_a(H_t;t,s,X_s)=\frac{1}{\Gamma(\frac{\beta+1}{2})}
\int_0^t\frac{G_a(H_t;t,s,X_s)}{(t-s)^{(1-\beta)/2}}dW_s^a.$$

One application of fractional stochastic equation of the form
(\ref{eq2}) is in finance. Fractional Black-Scholes market is
described in terms of the bank account and a stock. The price at a
time $t$ is given by the following formula
\begin{equation}
A_t=exp\Big(\int_0^t r(s)ds\Big),
\end{equation}where $r(s)\geq0,$ $s\in[0,t],$ represents the
interest rate. The price can be expressed using a fractional
Volterra-type equation:
\begin{equation}
X_t=X_0+\frac{1}{\Gamma(\beta)}\int_0^t\frac{\mu(s)X_s}{(t-s)^{1-\beta}}ds+
\frac{1}{\Gamma(\frac{\beta+1}{2})}
\int_0^t\frac{\sigma(s)X_s}{(t-s)^{(1-\beta)/2}}dW_s,
\end{equation}where $\mu,\sigma\geq0$ are continuous functions on
$[0,T].$

In this paper, we restrict our attention to stochastic fractional
Hamiltonian systems characterized by Wiener processes and assume
that the space of admissible curves in configuration space is of
class $\mathrm{C}^1.$ Random effects appear in the balance of
momentum equations, as white noise, that is why we may consider
randomly perturbed mechanical systems. It should be mentioned that
the ideas in this paper can be readily extended to stochastic
Hamiltonian systems \cite{milstein} driven by more general
semimartingales, but for the sake of clarity we restrict to Wiener
processes. Within this context, the results of the paper are as
follows:

\par 1. The paper proves almost surely that a  curve satisfies
stochastic fractional  HP equations if and only if it extremizes a
stochastic action. This theorem is the main result of the paper;

\par 2. Fractional HP equations are described using fractional Riemann-Liouville integral
and fractional It\^o integral;

\par 3. Langevian type stochastic fractional equations are
obtained, in the case of a hyperregular Lagrange function.

The paper is organized as follows. In Section 2 we present some
sufficient conditions for existence, uniqueness and almost sure
differentiability of stochastic flows on manifolds. In Section 3,
we extend the fractional Hamilton-Pontryagin (HP) principle to the
stochastic setting to prove that a class of mechanical systems
with multiplicative noise appearing as  forces and torques possess
a variational structure. In Section 4, for a hyperregular Lagrange
function, we get the stochastic fractional Hamiltonian equations
that lead to Langevin fractional equations. For a fractional
Lagrange function, defined on $\mathbb{R}^2,$ the corresponding
fractional Langevin equations are simulated. The mechanical system
could evolve on a nonlinear configuration space and involve
holonomic constraints or nonconservative effects in the drift. The
fractional Hamiltonian and the Lagrangian description are joined
together to get the fractional HP system.

\section{Stochastic flows on manifolds}

Some standard results on flows of SDE on manifolds are reviewed
here to the reader's convenience. For more detailed exposition,
the reader is referred to the textbooks  such as \cite{emery} or
\cite{ikeda}. This section parallels the treatment of
deterministic flows on manifolds found in \cite{abraham}.

Let $M$ be  a manifold modelled on a Banch space $E.$ Recall that
a vector field on the manifold $M$ is a section of the tangent
bundle $TM$ on $M.$ The set of all $\mathrm{C}^k$ vector fields on
$M$ is denoted by $\mathcal{X}^k(M).$

A stochastic dynamical system consists of a base flow on the
probability space which propagates the noise, and a stochastic
flow on $M$ which depends on the noise.

A stochastic dynamical system consists of a base flow on the
probability space $(\Omega,\mathcal{F},P)$ and a stochastic flow
on a manifold $M.$ The \emph{base flow} is a $P-$preserving map
$\theta:\mathbb{R}\times \Omega\rightarrow\Omega$ which satisfies:
\begin{enumerate}
    \item $\theta_0=id_{\Omega}:\Omega\rightarrow\Omega$ is the
    identity on $\Omega;$
    \item for all $s,t\in\mathbb{R},$ the group property,$\theta_s\circ\theta_t=
    \theta_{t+s}.$
\end{enumerate}

Given times $0\leq r\leq s\leq t,$ the "stochastic flow" on $M$ is
a map $\varphi_{t,s}:\Omega\times M\rightarrow M$ such that
\begin{enumerate}
    \item for almost all $\omega\in \Omega,$ the map $(s,t,\omega,x)\mapsto
    \varphi_{t,s}(\omega)x$ is continuous in $s,t$ and $x;$
    \item $\varphi_{s,s}(\omega)=id_{M}:M\rightarrow M$ is the
    identity map on $M,$ for all $s\in \mathbb{R};$
    \item $\varphi$ satisfies the cocycle property
       $$\varphi_{t,s}(\theta_s(\omega))\circ\varphi_{s,r}(\omega)=\varphi_{t,r}(\omega).$$
\end{enumerate}

This paper is concern with stochastic dynamical systems that come
from stochastic laws of motion, i.e. ones whose stochastic flows
define solutions of stochastic differential equations. Consider a
manifold $M,$ modelled on a Banach space $E$ and vector fields
$X_0, X_a\in \mathcal{X}^k(M), \, a=1,...,m.$ Let
$(W^a(t,\omega),\mathcal{F}_\tau), \, a=1,...,m,$ be independent
Wiener processes for $0\leq t \leq T.$ In terms of these objects,
the Stratonovich stochastic differential equations, that the paper
considers, takes the form:
\begin{equation}\label{1}
dx=X_0(x)dt+X_a(x)\circ dW^a, \, x(0)=x_0.
\end{equation}

A \emph{Stratonovich integral curve} of (\ref{1}) is a
$\mathrm{C}^0-$map, $c(\cdot,\omega):[0,T]\rightarrow M$ which
satisfies
$$c(t,\omega)=x_0+\int_0^t X_0(c(s,\omega))ds+\int_0^t X_a(c(s,\omega))\circ dW^a(s,\omega)ds,$$
for all $t\in[0,T].$

Let $c$ be a Stratonovich integral curve of (\ref{1}).
\emph{Pathwise uniqueness} of $c$ means that if
$\bar{c}:I\rightarrow M$ is also a solution of (\ref{1}) on the
same filtered probability space, with the same Brownian motion and
initial random variables, then,
$$P(c(t,\omega)=\bar{c}(t,\omega), \quad \forall \, t\in[0,T])=1. $$

For the rest of the paper the explicit dependence of stochastic
maps on the point $\omega\in\Omega$ will usually be suppressed.
With these definitions, one can state the following key, but
standard theorem (\cite{elwo}, \cite{emery}, \cite{ikeda}).

\begin{theorem}
(Existence, uniqueness and smoothness)\\
Let $M$ be a manifold with the model space $E.$ Suppose that
$X_0,X_a\in \mathcal{X}^k(M), \, a=1,...,m$ and $k\geq 1$ are
uniformly Lipschitz and measurable with respect to $x\in M.$ Let
$I=[0,T].$ Then the following statements hold.
\begin{enumerate}
    \item For each $u\in M,$ there is almost surely a
    $\mathrm{C}^0-$curve, $c:I\rightarrow M,$ such that $c(0)=u$
    and $c$ satisfies \emph{(\ref{1})}, for all $t.$ This curve $c:I\rightarrow M$
    is called a \emph{maximal solution};
    \item $c$ is pathwise unique;
    \item There is almost surely a mapping $F:I\times M\rightarrow
    M,$such that the curve $c_u:I\rightarrow M,$ defined by
    $c_u(t)=F_t(u),$is a curve satisfying \emph{(\ref{1})}, for
    all $t\in I.$ Moreover, almost surely $F$ is $\mathrm{C}^k$ in
    $u$ and $\mathrm{C}^0$ in $t.$
\end{enumerate}\hfill $\Box$
\end{theorem}

\section{Stochastic fractional HP mechanics}

In this section, a fractional variational principle is introduced
for a class of stochastic fractional Hamiltonian systems on
manifolds. The stochastic fractional action is a sum of classical
fractional action and several stochastic integrals.

Let $Q$ be an $n-$dimensional manifold, $\gamma_a:Q\rightarrow
\mathbb{R}, \, a=1,...,m,$ a deterministic function, and the
Lagrangian $L:TQ\rightarrow \mathbb{R}.$ Let $(\Omega,
\mathcal{F}, P)$ be a probability space and the interval
$[a,b]\subset \mathbb{R}.$ Let $\{W^a(t), \mathcal{F}_t\}_{t\in
[a,b]},$ for $a=1,...,m,$ where $\{W^a\}_{a=1,...,m}$ are
independent, real-valued Wiener processes and $\mathcal{F}_t$ is
the filtration generated by these Wiener Processes.

The stochastic HP fractional action is defined by
$\mathcal{A}^{\alpha,\beta}:\Omega\times
\mathcal{C}(PQ)\rightarrow\mathbb{R}$:
\begin{equation}\label{2}
\begin{array}{ll}
\mathcal{A}^{\alpha,\beta}(q,v,p,t)=\frac{1}{\Gamma(\alpha)}
\int_a^b(L(q(s),v(s))(t-s)^{\alpha-1}
+(p(s),\dot{q}(s)-v(s))(t-s)^{\alpha-1})ds+\\
\quad \\
\quad \quad \quad \quad \quad \quad
+\frac{1}{\Gamma(\beta)}\int_a^b\gamma_a(q(s))
(t-s)^{\beta-1}\circ dW^a(ds),
\end{array}
\end{equation}
where $0<\alpha\leq 1,$ $0<\beta\leq 1,$ $\Gamma(\alpha), \,
\Gamma(\beta)$ are the Euler gamma functions and $PQ=TQ\oplus
T^*Q,$ and
\begin{equation}\label{3}
\mathcal{C}(PQ)=\{(q,v,p)\in \mathrm{C}^0([a,b],PQ)|q\in
\mathrm{C}^1([a,b],Q), \, q(a)=q_a, \, q(b)=q_b\}.
\end{equation}

The first integral in (\ref{2}) is a Riemann integral, and the
second one is an It\^o integral.

The HP path space is a smooth infinite-dimensional manifold. One
can show that its tangent space in
$c=(q,v,p)\in\mathcal{C}([a,b],q_1,q_2)$ consists of maps
$w=(q,v,p,\delta q,\delta v,\delta p)\in
\mathrm{C}^0([a,b],T(PQ)),$ such that $\delta q_a=\delta q_b=0,$
and $q, \, \delta q$ are of class $\mathrm{C}^1.$ Let us denote by
$(q,v,p)(\cdot,\varepsilon)\in \mathcal{\mathcal{C}}(PQ)$ the
one-parameter family of curves in $\mathcal{C},$ that is
differentiable with respect to $\varepsilon.$ Define the
differential of $\mathcal{A}^{\alpha,\beta}$ as
\begin{equation}\label{4}
d\mathcal{A}^{\alpha,\beta}(\delta q,\delta v,\delta
p)=\frac{\partial}{\partial
\varepsilon}\mathcal{A}^{\alpha,\beta}(\omega,q(s,\varepsilon),
v(s,\varepsilon),p(s,\varepsilon))\Big|_{\varepsilon=0},
\end{equation}
where
\begin{equation}\label{5}
\delta q(s)=\frac{\partial}{\partial \varepsilon}
q(s,\varepsilon)\Big|_{\varepsilon=0}, \, \delta q(a)=\delta
q(b)=0, \, \delta v(s)\frac{\partial}{\partial \varepsilon}
v(s,\varepsilon)\Big|_{\varepsilon=0}, \,
dp(s)=\frac{\partial}{\partial \varepsilon}
p(s,\varepsilon)\Big|_{\varepsilon=0}.
\end{equation}

In terms of this differential, one can state the following
critical point condition for the action
$\mathcal{A}^{\alpha,\beta}.$

\begin{theorem}
(Stochastic fractional variational principle of HP)\\
Let $L:TQ\rightarrow \mathbb{R}$ be a Lagrangian on $TQ$ of class
$\mathrm{C}^2,$ with respect to $q$ and $v,$ and with globally
Lipschitz first derivative with respect to $q$ and $v.$ Let
$\gamma_a:Q\rightarrow \mathbb{R}$ be of class $\mathrm{C}^2$ and
with globally Lipschitz first derivatives, for $a=1,...,m.$ Then,
almost surely, a curve $c=(q,v,s)\in\mathcal{C}(PQ)$ satisties the
stochastic fractional equations
\begin{equation}\label{6}
\begin{array}{ll}
dq(s)=v(s)ds,\\
\quad \\
dp(s)=\Big(\frac{\partial L}{\partial q}(q(s),v(s))+\frac{\partial
L}{\partial
v}(q(s),v(s))\frac{\alpha-1}{t-s}\Big)ds+\frac{\partial
\gamma_a(q(s))}{\partial
q}\frac{\Gamma(\alpha)}{\Gamma(\beta)}(t-s)^{\beta-\alpha}\circ
dW^a(s),\\
\quad \\
p(s)=\frac{\partial L}{\partial v}(q(s),v(s)), \, s\in[a,b],
\end{array}
\end{equation}if and only if it is a critical point of the function
$\mathcal{A}^{\alpha,\beta}:\Omega\times
\mathcal{C}(PQ)\rightarrow \mathbb{R},$ i.e.
$d\mathcal{A}^{\alpha,\beta}(c)=0.$
\end{theorem}
\textbf{Proof:} The proof results by applying the method from \cite{bou} and \cite{el}.\hfill $\Box$\\

Observe that by the It\^o-Stratonovich conversion formula, the
It\^o modification to the drift is equal to 0, and hence (\ref{6})
can be written in the It\^o form as
\begin{equation}\label{7}
\begin{array}{ll}
dq=vds,\\
\quad \\
dp=\Big(\frac{\partial L}{\partial q}(q,v)+\frac{\partial
L}{\partial v}(q,v)\frac{\alpha-1}{t-s}\Big)ds+\frac{\partial
\gamma_a(q)}{\partial
q}\frac{\Gamma(\alpha)}{\Gamma(\beta)}(t-s)^{\beta-\alpha}
dW^a(s),\\
\quad \\
p=\frac{\partial L}{\partial v}(q,v), \, s\in[a,b].
\end{array}
\end{equation}

In what follows, structure-preserving properties  of the flow map,
defined by maximal solution of the equations over $[a,b]$ will be
investigated. First, observe that because of smoothness conditions
assumed in Theorem 2, a solution almost surely exists and it is
pathwise unique on $[a,b],$ by the result from Section 2. When
$\gamma_a$ is constant for, $a=1,...,m$ and $\alpha=1,$ the reader
is referred to \cite{yoshi}, for deterministic treatments of
symplecticity, momentum map preservation and holonomically
constrained mechanical systems.

If $\beta=\alpha,$ the equation (\ref{7}) is given by:
\begin{equation}\label{8}
\begin{array}{ll}
dq=vds,\\
\quad \\
dp=\Big(\frac{\partial L}{\partial q}(q,v)+\frac{\partial
L}{\partial v}(q,v)\frac{\alpha-1}{t-s}\Big)ds+\frac{\partial
\gamma_a(q)}{\partial q}
dW^a(s),\\
\quad \\
p=\frac{\partial L}{\partial v}(q,v).
\end{array}
\end{equation}

If $\gamma_a$ is constant, for $a=1,...,m,$ from (\ref{7}),
results:
\begin{equation}\label{9}
\begin{array}{ll}
dq=vds,\\
\quad \\
dp=\Big(\frac{\partial L}{\partial
q}(q,v)+\frac{\partial L}{\partial v}(q,v)\frac{\alpha-1}{t-s}\Big)ds,\\
\quad \\
p=\frac{\partial L}{\partial v}(q,v),
\end{array}
\end{equation}
and they represent the Euler-Lagrange fractional equations
(\cite{el}, \cite{federico}).

\section{Stochastic fractional equation for the Lagrangian hyperregular}

Let $L:TQ\rightarrow \mathbb{R}$ be a Lagrangian on $TQ$
hyperregular, that means $det\Big(\frac{\partial^2 L}{\partial v^i
\partial v^j}\Big)\neq0.$ From (\ref{6}) results the following
propositions:

\begin{pr}
(Stochastic fractional Hamiltonian equations)\\
The equations \emph{(\ref{7})} are equivalent with the equations:
\begin{equation}\label{14}
\begin{array}{ll}
dq^i=\frac{\partial H}{\partial p_i}ds,\\
\quad \\
dp_i=\Big(-\frac{\partial H}{\partial
q^i}+\frac{\alpha-1}{t-s}p_i\Big)ds+\frac{\Gamma(\alpha)}{\Gamma(\beta)}\frac{\partial
\gamma_a(q)}{\partial q^i}(t-s)^{\beta-\alpha}dW^a(s), \,
i=1,...,n,
\end{array}
\end{equation}where $H=p_iq^i-L(q,v).$ \hfill $\Box$
\end{pr}

\begin{pr}
If $L=\frac{1}{2}g_{ij}v^iv^j,$ where $g_{ij}$ are the components
of a metric on the manifold $Q,$ equations \emph{(\ref{7})} take
the form:
\begin{equation}\label{15}
\begin{array}{ll}
dq^i=v^ids,\\
\quad\\
dv^i=-\Big(\Gamma_{jk}^i
v^jv^k+\frac{\alpha-1}{t-s}v^i\Big)ds+\frac{\Gamma(\beta)}{\Gamma(\alpha)}g^{ij}
\frac{\partial \gamma_a(q)}{\partial q^j}(t-s)^{\beta-1}dW^a(s),
\, i,j=1,...,n,
\end{array}
\end{equation}
where $\Gamma_{jk}^i$ are Cristofel coefficients associated to the
considered metric. Equations \emph{(\ref{14})} become:
\begin{equation}\label{16}
\begin{array}{ll}
dq^i=v^ids,\\
\quad \\
dp_i=\Big(\frac{1}{ 2}\frac{\partial g_{kl}}{\partial
q^i}p^kp^l+\frac{\alpha-1}{t-s}p_i\Big)ds+\frac{\Gamma(\alpha)}{\Gamma(\beta)}\frac{\partial
\gamma_a(s)}{\partial q^i}(t-s)^{\beta-\alpha}dW^a(s), \,
i,j,k=1,...,n,
\end{array}
\end{equation}\hfill $\Box$
\end{pr}
Equations (\ref{14}) represent fractional Langevin equations.
Equations (\ref{16}) can be used for fractional motion of
relativistic particles with noise.

\begin{pr}
\begin{description}
    \item[a)] If $Q=\mathbb{R},$  $H(p,q)=\frac{1}{2}p^2+U(q)$
    and $\gamma(q)=\cos(q),$ equations \emph{(\ref{14})} are given
    by:
    \begin{equation}\label{17}
    \begin{array}{ll}
dq=pds,\\
\quad\\
dp=\Big(-\frac{dU}{dq}+\frac{\alpha-1}{t-s}p\Big)ds-\frac{\Gamma(\alpha)}{\Gamma(\beta)}
(t-s)^{\beta-\alpha}\sin(q)dW(s);
    \end{array}
    \end{equation}
    \item[b)] If $U(q)=\cos(q),$ the Euler scheme for
    \emph{(\ref{17})}is:
    \begin{equation}\label{18}
    \begin{array}{ll}
x(n+1)=x(n)+hy(n),\\
y(n+1)=y(n)+h(\sin(x(n))+\frac{\alpha-1}{t-nh}y(n))-\frac{\Gamma(\alpha)}{\Gamma(\beta)}
(t-nh)^{\beta-\alpha}sin(x(n))G(n),\\
    \end{array}
    \end{equation}where $h=\frac{T}{N},$ $G(n)=W((n+1)h)-W(nh),$
    $n=0,...,N-1,$ $0<\alpha<1,\, 0<\beta<1,$ and $x(n)=q(nh),\, y(n)=p(nh).$
\end{description}
\end{pr}

For $\alpha=0.6, \, \beta=0.3, \, t=0.8$ and $h=0.0001,$ with
Maple 13,  the orbit $(n,p(nh))$ is represented in Figure 1, and
the orbit $(n,p(nh,\omega))$ is represented in Figure 2

\begin{center}\begin{tabular}{cc}
\epsfxsize=6cm \epsfysize=5cm
 \epsffile{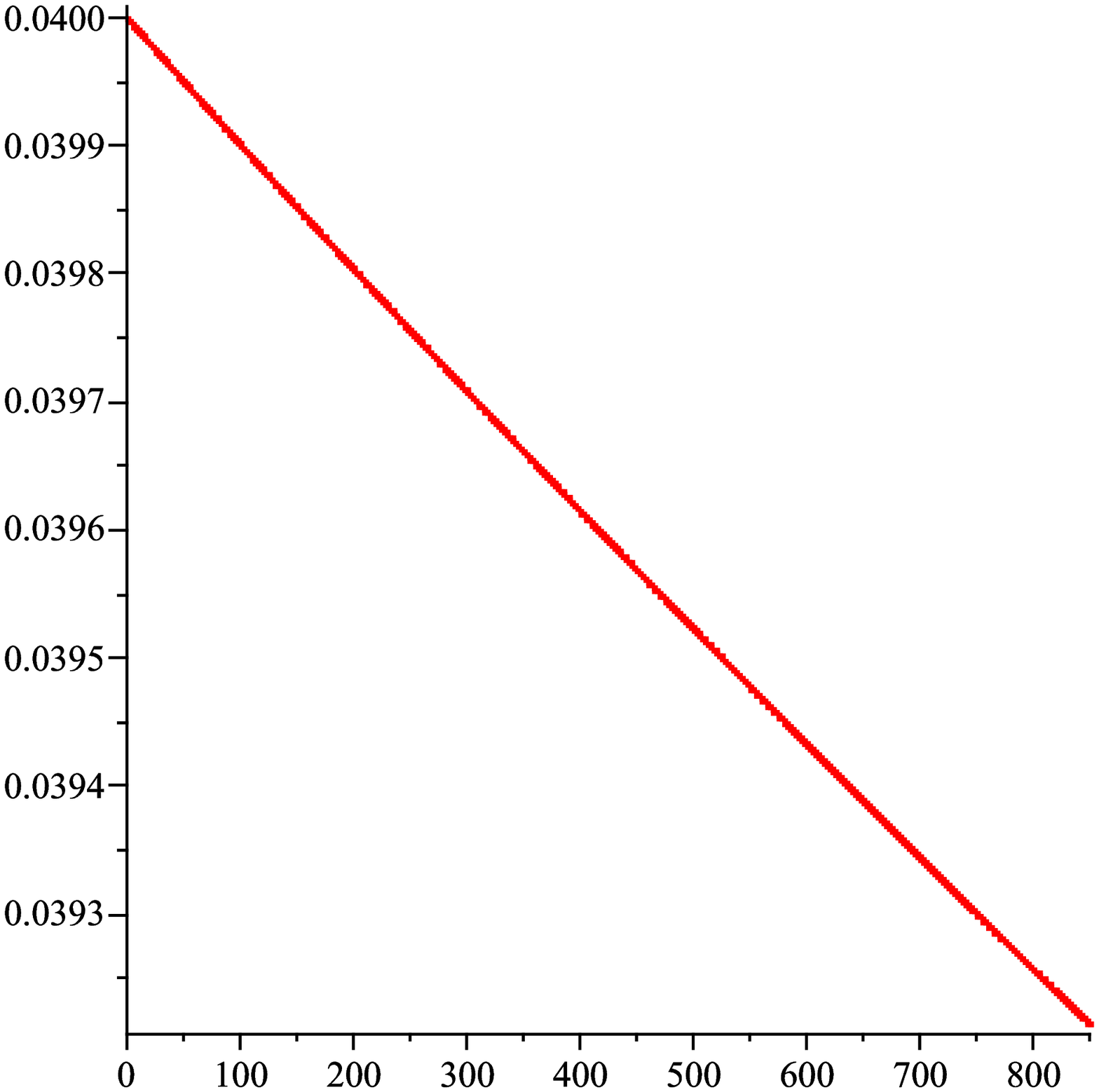}  &
\epsfxsize=6cm \epsfysize=5cm
\epsffile{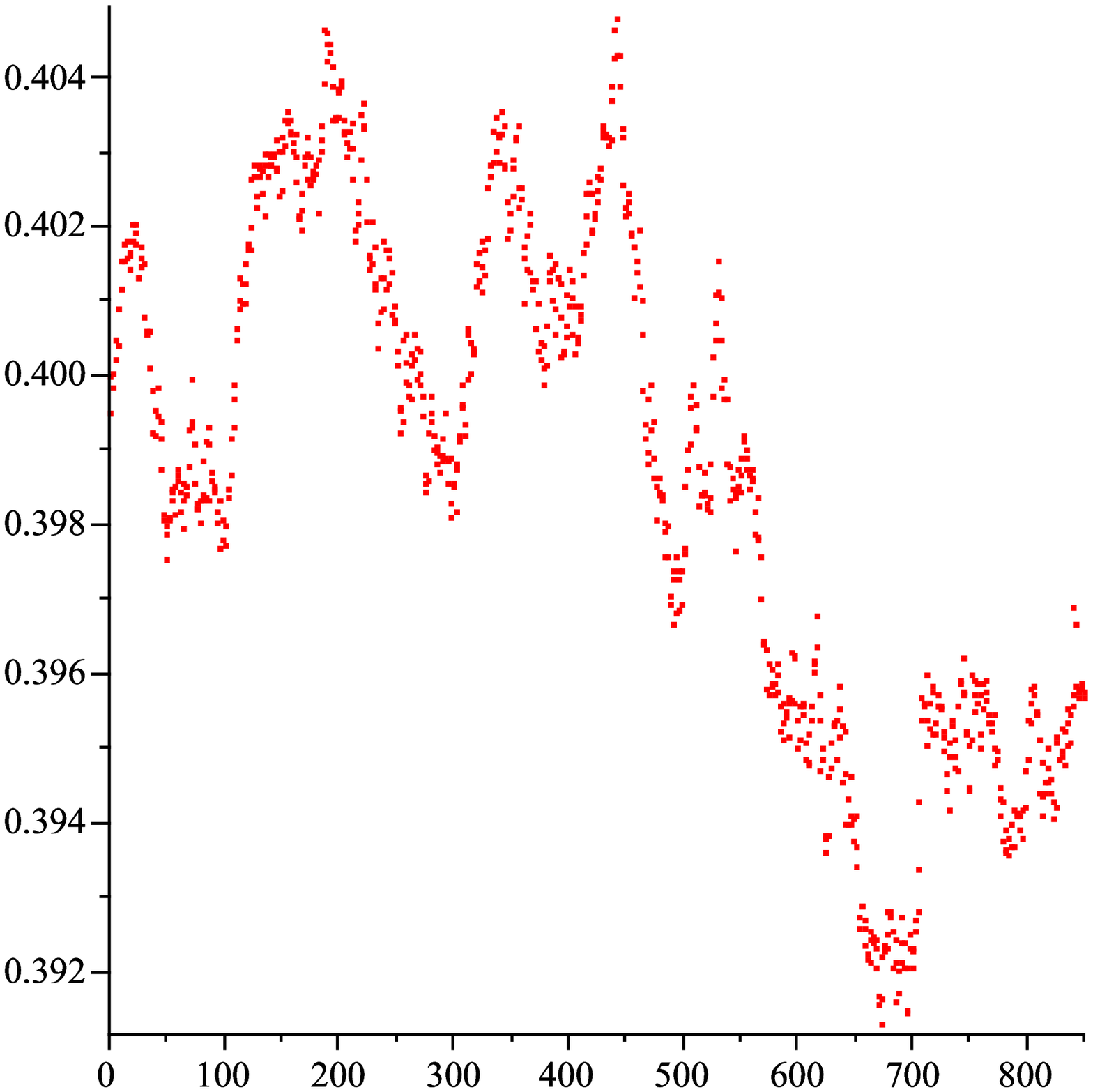} \\
Figure 1: the orbit $(n,p(nh))$&
Figure 2: the orbit $(n,p(n,\omega))$\\
    \end{tabular}
\end{center}

In Figure 3 the orbit $(q(nh),p(nh))$ is represented, and in
Figure 4 it is represented the orbit
$(q(nh,\omega),p(nh,\omega)).$

\begin{center}\begin{tabular}{cc}
\epsfxsize=6cm \epsfysize=5cm
 \epsffile{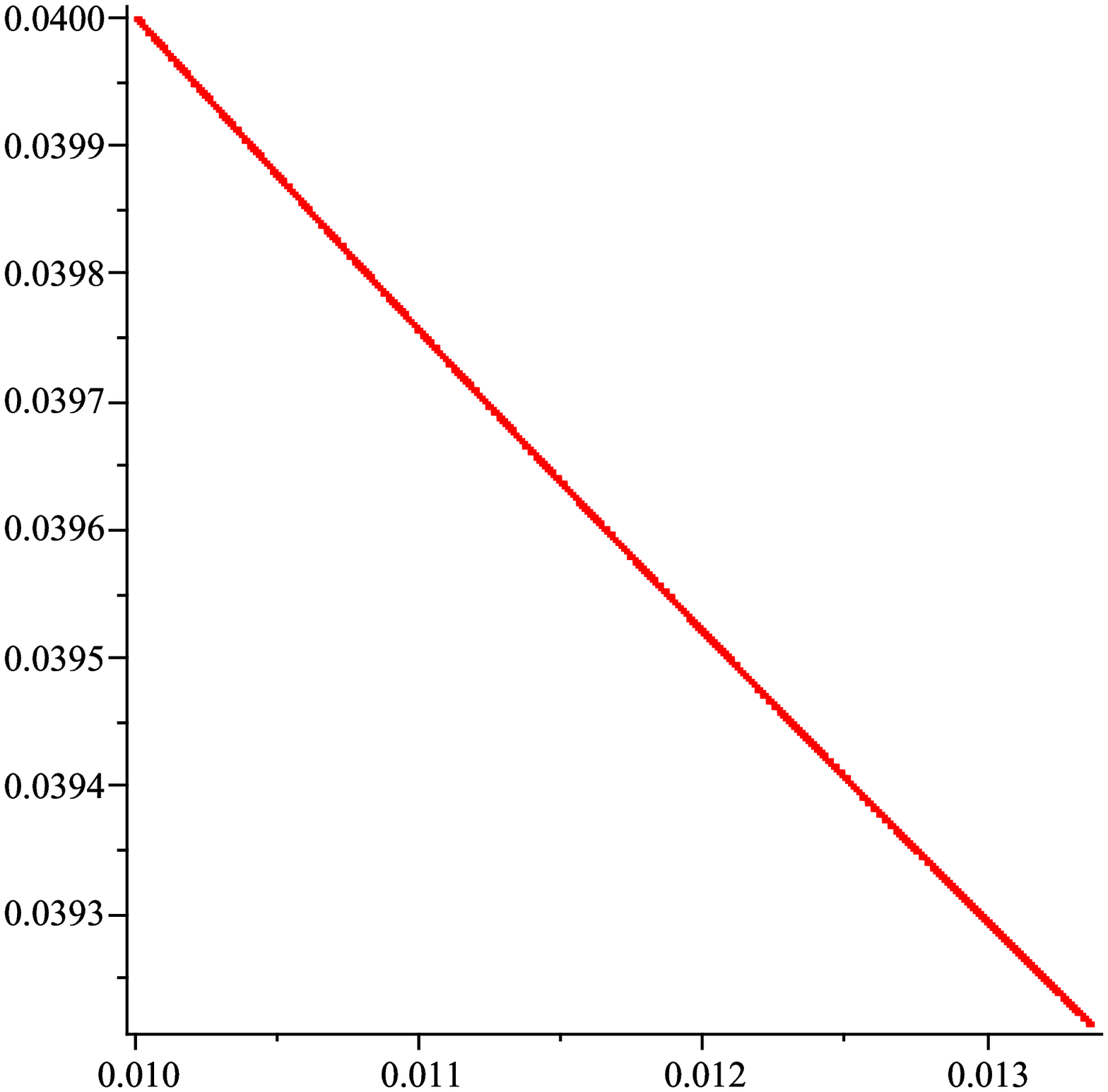}  &
\epsfxsize=6cm \epsfysize=5cm
\epsffile{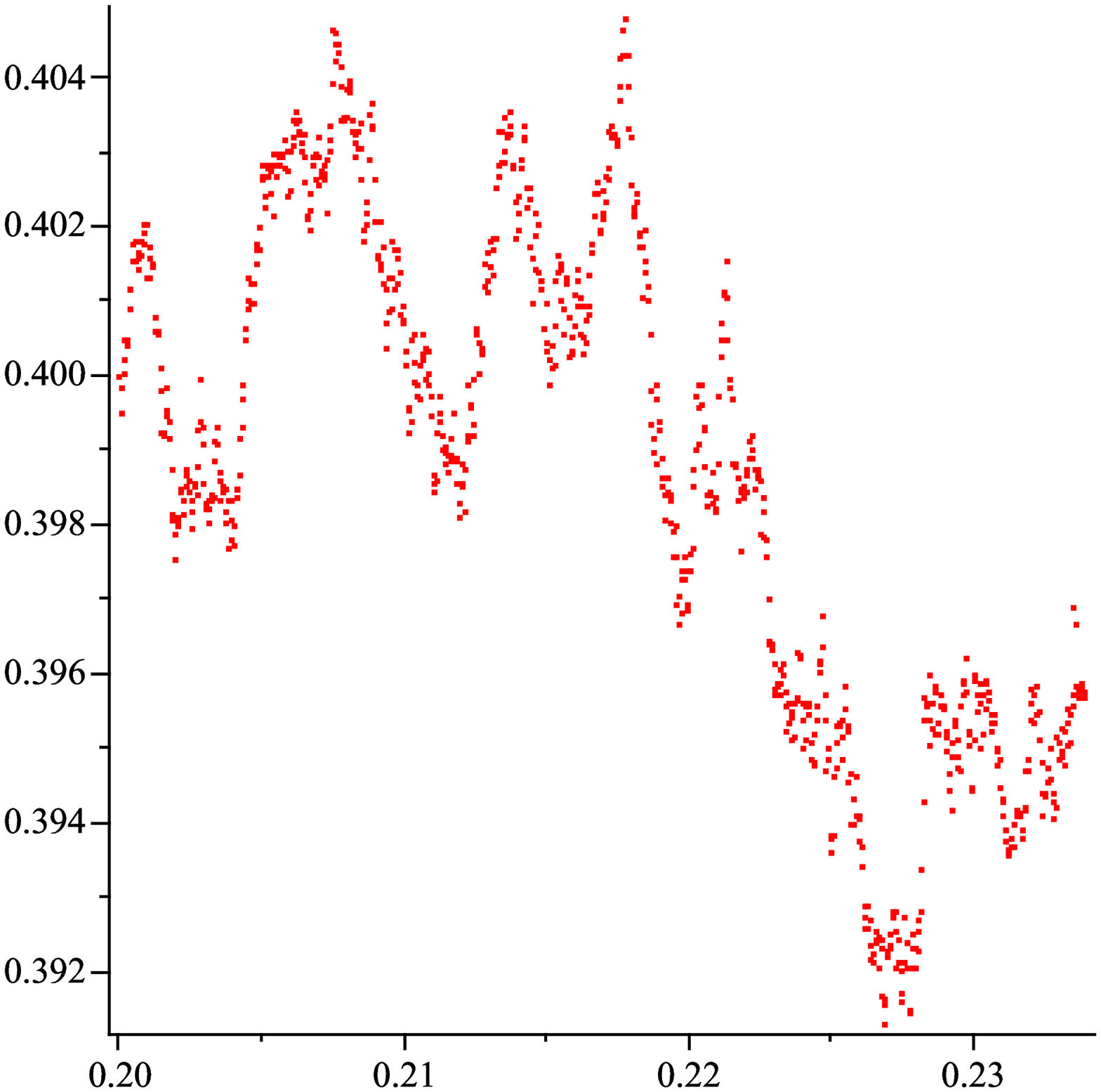} \\
Figure 3: the orbit $(q(nh),p(nh))$&
Figure 4: the orbit $(q(nh,\omega),p(nh,\omega))$\\
    \end{tabular}
\end{center}

For $\alpha=\beta=0.6,$ and $t=0.8,\, h=0.0001,$ the orbits
$(n,p(nh,\omega))$ and $(q(nh,\omega),p(nh,\omega))$ are
represented in Figure 5 and Figure 6.

\begin{center}\begin{tabular}{cc}
\epsfxsize=6cm \epsfysize=5cm
 \epsffile{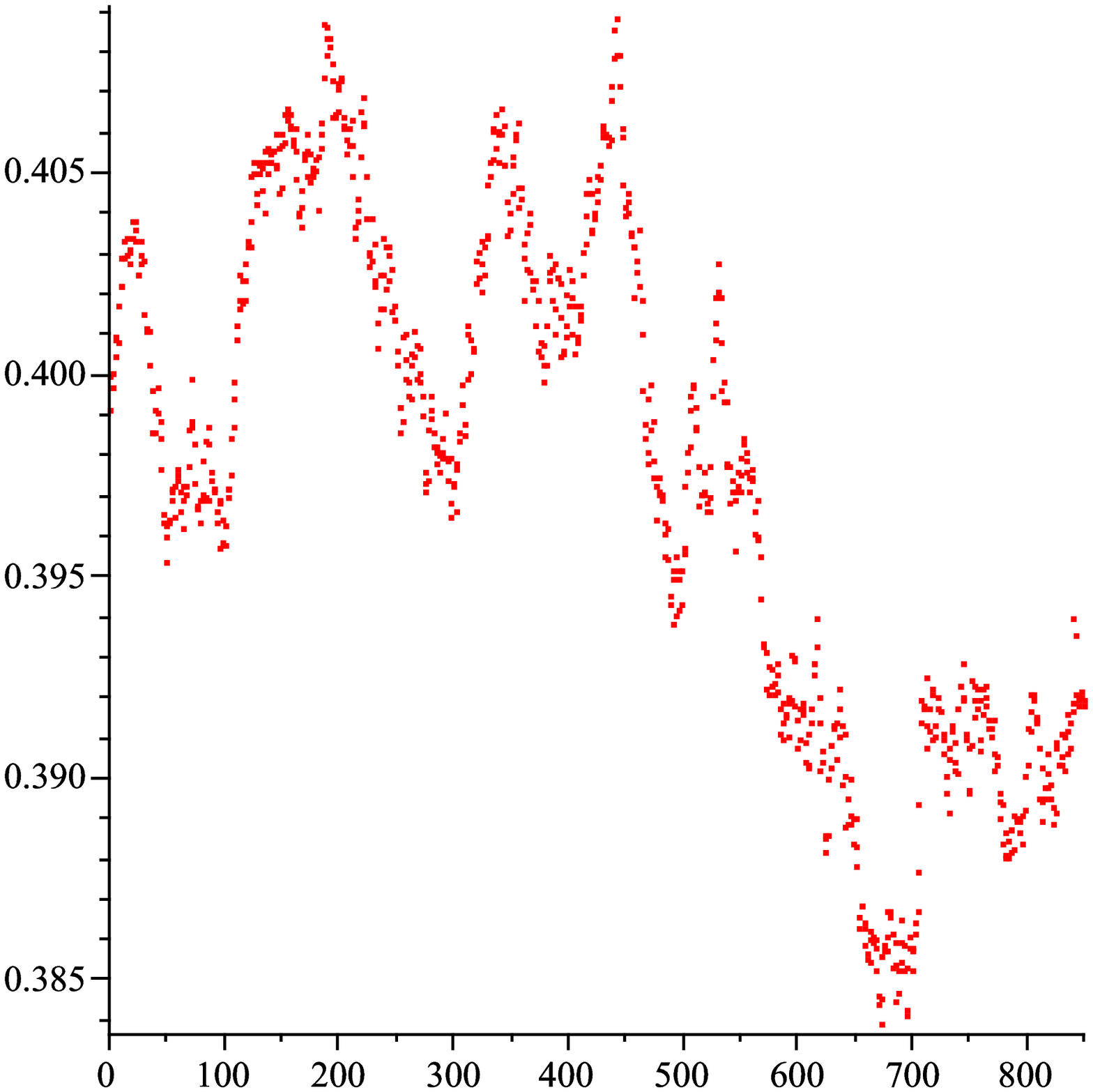}  &
\epsfxsize=6cm \epsfysize=5cm
\epsffile{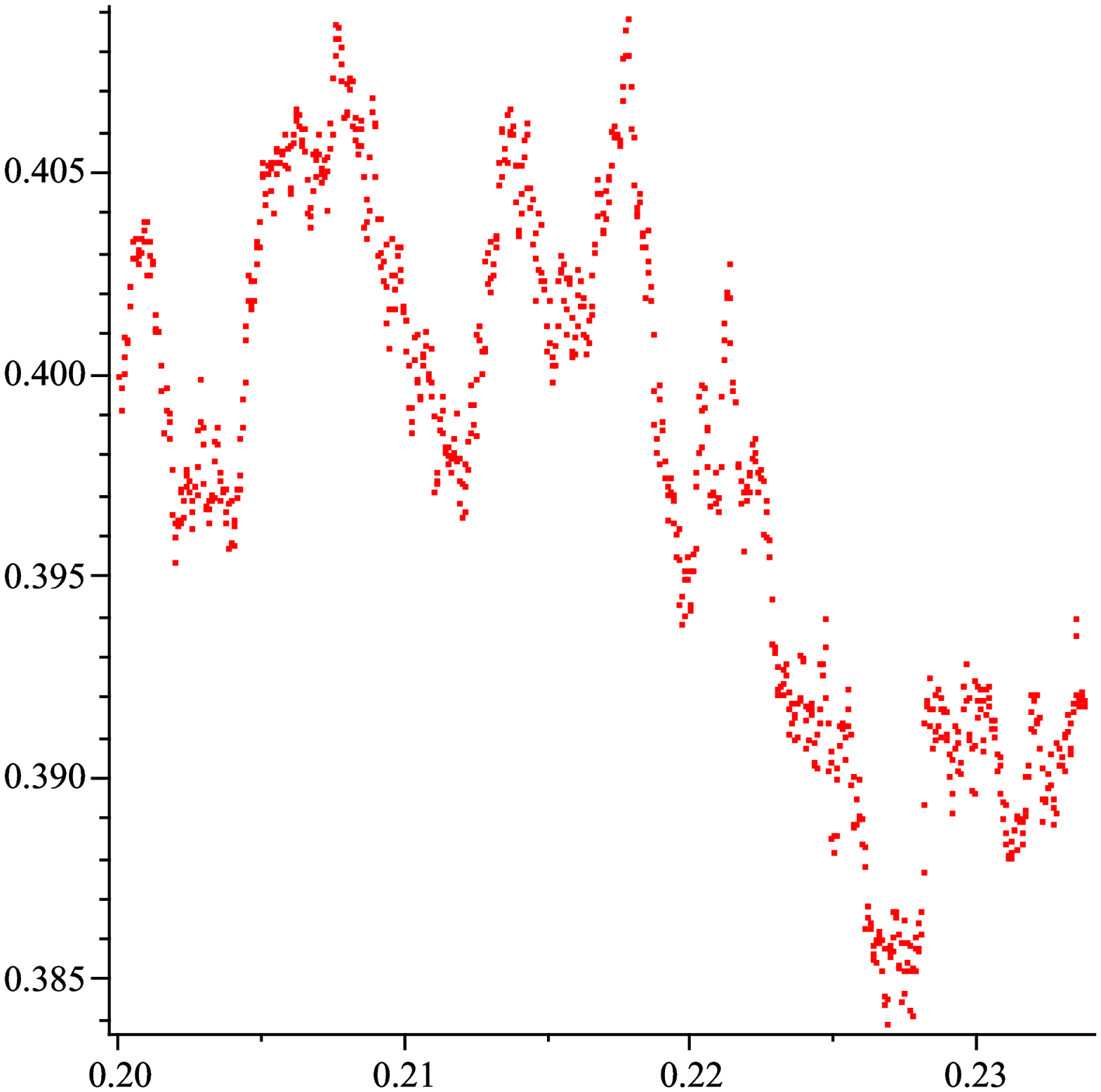} \\
Figure 5: the orbit $(n,p(nh,\omega))$&
Figure 6: the orbit $(q(nh,\omega),p(nh,\omega))$\\
    \end{tabular}
\end{center}

\section*{Conclusions}

In this paper it was described stochastic fractional HP principle,
using classical stochastic HP principle \cite{bou} and fractional
principle (\cite{el}, \cite{federico}). Using a hyperregular
Lagrange function, Langevin-type fractional equations were
illustrated. We have done the numerical simulation for the case of
a Lagrange function defined on $\mathbb{R}^2.$

In the future work, we will consider other problems that deal with
stochastic fractional HP principle.


\begin{thebibliography}{99}


\bibitem{abraham} Abraham, R., Marsden, J.E., Ra\c{t}iu, T., \emph{Manifolds,
Tensors, Analysis, and Applications}, New York, Springer-Verlag,
2007.

\bibitem{agrawal1} Agrawal, O.P., \emph{ Formulation of Euler-Lagrange equations for fractional
variational problems}, J. Math. Anal. Appl. 272 (2002), no. 1,
368-379.

\bibitem{agrawal} Agrawal, O.P., \emph{Formulation of Euler-Lagrange equations
for fractional variational problems}, J. Math. Anal. Appl. 272
(2002), no. 1, 368-379.

\bibitem{baleanu} B\u{a}leanu, Avkar, T., \emph{Lagrangians with linear velocities within
Riemann-Liouville fractional derivatives}, Nuovo cimento 119,
(2004) 73-79.

\bibitem{bou} Bou-Rabee, N., \emph{Stochastic variational integrators},
 IMA Journal of Numerical Analysis Advance, 2008.


\bibitem{chis} Chi\c{s}, O., Despi, I., Opri\c{s}, D., \emph{Fractional
  equations on algebroids and fractional
algebroids},vol. New Trends in Nanotechnology and Fractional
Calculus Applications, Springer-Verlag, Berlin, Heidelberg, New
York, will apear.

\bibitem{chis2} Chi\c{s}, O., Opri\c{s}, D., \emph{Mathematical pendulum and its
variants}, arXiv:0905.4356v1[math.DS].

\bibitem{chis3} Chi\c{s}, O., Opri\c{s}, D., \emph{Mathematical analysis of stochastic
models for tumor-immune systems}, arXiv:0906.2794v1[math.DS] (sent
for publication).

\bibitem{el-b1} El-Borai, M.M., El-Said El-Nadi, O.L., Mostafa,
Ahmed, H.H., \emph{Volterra equations with fractional stochastic
integrals},Mathematical problems in Engineering, 5 (2004),
453-468.

\bibitem{el} El-Nabulsi, R.A.,  \emph{A fractional action-like variational
approach of some classical, quantum and geometrical dynamics},
Int. J. Appl. Math. 17 (2005), 299-317.

\bibitem{elwo} Elworty, K.D., \emph{Stochastic Differential Equations on
Manifolds}, Cambridge, UK: Cambridge University Press, 1982.

\bibitem{emery} Emery, M., \emph{Stochastic Calculus in Manifolds},
Berlin, Springer-Verlag, 1989.

\bibitem{federico} Frederico, G.S.F., Torres, D.F.M,  \emph{A formulation of Noether's
theorem for fractional problems of the calculus of variations}, J.
Math. Anal. Appl. 334 (2007), no. 2, 834-846.

\bibitem{gore} Gorenflo, R., Mainardi, F., \emph{Fractional calculus and stable
probability distributions}, Arch Mech 1995;50(3):377-88.

\bibitem{ikeda} Ikeda, N., Watabe, S., \emph{Stochastic Differential
Equations and Diffusion Processes}, Amsterdam, North-Holland,
1989.

\bibitem{klimek}  Klimek, M.,  \emph{Lagrangean and Hamiltonian fractional sequential
mechanics}, Czechoslovak J. Phys. 52 (2002), no. 11, 1247-1253.

\bibitem{lazaro} Lazaro-Cami, J.A., Ortega, J.P., \emph{Reduction and
reconstruction of symmetric stochastic differential equations},
Rep. Math. Phys., a2007, in press.

\bibitem{lazaro1} Lazaro-Cami, J.A., Ortega, J.P., \emph{Stochastic
Hamiltonian dynamical systems}, Rep. Math. Phys, b2007, in press.

\bibitem{milstein} Milstein, G.N., Repin, YU.M., Tretyakov, M.V.,
\emph{Symplectic methods for Hamiltonian systems with additive
noise}, SIAM J. Numer. Anal., 39 (2002), 1-9.

\bibitem{pardoux}  Pardoux, E., Protter, P., \emph{A two-sided stochastic integral and
its calculus}, Probab. Theory Related Fields 76 (1987), no. 1,
15-49.

\bibitem{podl} Podlubny, I., \emph{Fractional Differential Equations}, Acad. Press,
San Diego, 1999.

\bibitem{tarasov} Tarasov, V.E., \emph{Fractional variations for dynamical systems:
Hamilton and Lagrange Approaches}, Journal of Physics A 39, No.26
(2006), 8409-8425.

\bibitem{udriste} Udri\c{s}te, C., Opri\c{s}, D.,\emph{ Multi-time Euler-Lagrange-Hamilton
theory}, WSEAS Transactions on Mathematics Issue 1 volume 7
(2008), 19-30.

\bibitem{yoshi} Yoshimura, H., Marsden, J.E.,\emph{ Dirac structures and
Lagrangian mechanics part I: implicit Lagrangian systems}, J.
Geom. Phys., 57 (2006), 133-156.

\end{thebibliography}
\end{document}